\documentclass[12pt,reqno]{amsart}

\usepackage{amscd}
\usepackage{amssymb,latexsym}
\usepackage{epsfig}
\usepackage{float}
\usepackage{url}
\usepackage{graphicx,xcolor}
\usepackage{tikz}
\renewcommand{\nocite}[1]{{#1}}

\newcommand{\N}{{\mathbb N}}
\newcommand{\Z}{{\mathbb Z}}

\newcommand{\R}{{\mathbb R}}

\newcommand{\NN}{{\mathcal N}}

\newcommand{\www}{\widetilde}

\DeclareMathOperator{\Der}{Der}

\DeclareMathOperator{\Fix}{Fix}

\DeclareMathOperator{\id}{id}

\DeclareMathOperator{\sign}{sign}

\makeatletter 
\@namedef{subjclassname@2020}{\textup{2020} Mathematics Subject Classification}
\makeatother

\begin{document}

\theoremstyle{plain}
\newtheorem{lemma}{Lemma}[section]
\newtheorem{definition/lemma}[lemma]{Definition/Lemma}
\newtheorem{theorem}[lemma]{Theorem}
\newtheorem{proposition}[lemma]{Proposition}
\newtheorem{conjecture}[lemma]{Conjecture}
\newtheorem{conjectures}[lemma]{Conjectures}
\newtheorem{corollary}[lemma]{Corollary}

\theoremstyle{definition}
\newtheorem{definition}[lemma]{Definition}
\newtheorem{withouttitle}[lemma]{}
\newtheorem{remark}[lemma]{Remark}
\newtheorem{remarks}[lemma]{Remarks}
\newtheorem{example}[lemma]{Example}
\newtheorem{examples}[lemma]{Examples}
\newtheorem{notations}[lemma]{Notations}

\title[Generic games with two strategies for each player]
{Maximal number of mixed Nash equilibria in generic games
where each player has two pure strategies} 

\author{Claus Hertling and Matija Vuji\'c}

\address{Claus Hertling\\
Lehrstuhl f\"ur algebraische Geometrie, Universit\"at Mannheim,
B6 26, 68159 Mannheim, Germany}

\email{claus.hertling@uni-mannheim.de}

\address{Matija Vuji\'c\\
Z\"urich, Switzerland}


\date{September 25, 2026}

\thanks{We would like to thank a referee for careful
and detailed feedback and J. Kocutar and M. Verdonck for
an example which led us to modify a conjecture.}



\begin{abstract}
{\Small The number of Nash equilibria of the mixed extension of
a generic finite game in normal form is finite and odd.
This raises the question how large the number can be,
depending on the number of players and the numbers of their
pure strategies. Here we present a lower bound for the maximal
possible number in the case of $n$-player games where each player 
has two pure strategies. It is surprisingly close to a known 
upper bound.}
\end{abstract}

\maketitle

\tableofcontents

\setcounter{section}{0}

\section{Main result and its context}\label{c1}
\setcounter{equation}{0}

The set $\NN$ of Nash equilibria of the mixed extension
of a finite game in normal form is non-empty.
This was proved by \nocite{Nash (1951)}, and it was one early
instance of the relevance of the notion of Nash equilibria.
A natural question is how the set $\NN$ might look like.
For generic games, it was shown by \nocite{Wilson (1971)},
\nocite{Rosenmuller (1971)} and \nocite{Harsanyi (1973)}
that $\NN$ is finite and odd. 
For arbitrary games, $\NN$ is a real 
semialgebraic set and can look arbitrarily complicated 
(\nocite{Datta (2003)}).

Here we focus on the case of generic games and ask about
the maximum of the numbers $|\NN|$ for generic games
with fixed numbers of players and of pure strategies. 
This is not of immediate use in applications in economics,
but it is useful in order to estimate the complexity 
of the problem of controlling $\NN$.

We consider finite games with a fixed number $n\geq 2$ of players
$i\in[n]:=\{1,...,n\}$ where each player has a fixed set
$S^i$ of $m_i\in\N(=\{1,2,3,..\})$ pure strategies.
The set $S=S^1\times ...\times S^m$ is the
set of all pure strategy combinations. 
Each player has a utility function $U^i:S\to\R$.
These combine to a utility map $U=(U^1,...,U^m):S\to\R^m$.
A game is generic if $U$ is generic. 
In the mixed extension $([n],X,V)$, the set $X^i$ of mixed
strategies of player $i$ is the simplex of probability
distributions over $S^i$, and $V:X\to\R^m$ is the multilinear
extension of $U$ to $X=X^1\times ...\times X^m$. We denote by 
$\mu(n;m_1,...,m_n)\in\N$ the maximal number of mixed
Nash equilibria within all games $([n],S,U)$ with fixed sets 
$[n]$ and $S$ and generic $U$.

In most cases this number is unknown. In the case of two players
an upper bound was given by \nocite{Keiding (1997)}, and 
in the case of two players and $m_1=m_2$ 
a lower bound was given by \nocite{von Stengel (1997)}
and \nocite{von Stengel (1999)}.

For more than two players, little is known. The best result
is due to \nocite{McKelvey and McLennan (1997)}.
They considered {\it totally mixed Nash equilibria (TMNE)},
which are Nash equilibria in the interior of $X$,
so that any pure strategy has positive probability.
They found that the maximal number of TMNE for generic games 
$([n],S,U)$ as above is the number $E(n;m_1-1,...,m_n-1)$ of 
{\it block derangements}: It is the number of partitions of the set 
$\bigcup_{i\in[n]}(S^i-\{s^i_0\})$ (where $s^i_0$ is a fixed
strategy in $S^i$) into new sets $B^1,...,B^n$
with $|B^i|=m_i-1$ and $B^i\cap (S^i-\{s^i_0\})=\emptyset$. 
They and \nocite{Vidunas (2017)} studied these numbers.
Vidunas derived from this a rough upper bound 
for $\mu(n;m_1,...,m_n)$ in the general case. He used that
each Nash equilibrium is a TMNE of the restricted game where
only those pure strategies are considered which have positive
probability in the given Nash equilibrium.

In this paper we focus on the number $\mu(n;2,...,2)$,
so many players with two strategies each.
There the number of block derangements specializes
to the number 
\begin{eqnarray}\label{1.1}
!n:=E(n;1,...,1)
=|\{\sigma\in S_n\,|\, \textup{Fix}(\sigma)=\emptyset\}|
\end{eqnarray}
of {\it derangements} (permutations without fixed points). 
Vidunas' upper bound specializes to
(see Remark \ref{t2.2} (v) for $\leq$, see 
\nocite{Vidunas (2017)} for $\stackrel{!}{=}$) 
\begin{eqnarray}\label{1.2}
\mu(n;2,...,2)\leq V(n):=\sum_{l=0}^n 
\begin{pmatrix}n\\l\end{pmatrix}
\cdot 2^l\cdot !(n-l)\stackrel{!}{=}
\sum_{l=0}^n\frac{n!}{l!}.
\end{eqnarray}
Our main result is the lower bound
\begin{eqnarray}\label{1.3}
\frac{1}{2}(V(n)+!n)\leq \mu(n;2,...,2).
\end{eqnarray}
This is more than half of the upper bound $V(n)$,
so it is remarkably close to the upper bound. For example
\begin{eqnarray}\label{1.4}
\begin{array}{c|c|c|c|c|c}
n & 2 & 3 & 4 & 5 \\ \hline 
V(n) & 5 & 16 & 65 & 326 \\ \hline
\frac{1}{2}(V(n)+!n) & 
3 & 9 & 37 & 187
\end{array}
\end{eqnarray}
We even conjecture equality in \eqref{1.3} (Conjecture 
\ref{t4.5}).

Section \ref{c2} sets the notations for finite games 
where each player has two pure strategies and their mixed extensions.
Section \ref{c3} introduces a special subfamily of non-generic games,
and it gives general results on their Nash equilibria 
and {\it equilibrium candidates}.
It sets the stage for Section \ref{c4}. 
Section \ref{c4} proves that this subfamily contains games
with $\frac{1}{2}(V(n)+!n)$ Nash equilibria.
The art is to control simultaneously Nash equilibria with
all different possible supports. This leads to interesting
combinatorics.

\section{Two-action games and their mixed extensions}\label{c2}
\setcounter{equation}{0}

In this section, a finite game in normal form where each player
has two pure strategies, its mixed
extension, best reply maps, Nash equilibria and equilibrium
candidates are introduced formally.  
As usual, for two nonempty sets $A$ and $B$, $A^B$ denotes 
the set $\{\textup{maps from }B\textup{ to }A\}$.

\begin{definition}\label{t2.1}
(a) As in the introduction $n\in\Z_{\geq 2}$ and $[n]:=\{1,2,...,n\}$.
A {\it two-action game} is a game in normal form
$([n],S,U)$ with the following properties:
$[n]$ is the set of players. Player $i$ has the set $S^i=\{0;1\}$
of two pure strategies. $S=S^1\times ...\times S^n=\{0;1\}^{[n]}$ 
is the set of pure strategy combinations.
$U^i:S\to\R$ is the utility function of player $i$, and 
$U=U^1\times ...\times U^n:S\to\R^n$.

\medskip
(b) $([n],X,V)$ denotes the {\it mixed extension}
of the two-action game in part (a). We call also $([n],X,V)$ 
a {\it two-action game}. Here 
\begin{eqnarray*}
X^i&:=&[0,1]\quad\textup{ for }i\in[n]
\quad\textup{ with coordinate }x_i,\\
X&:=&X^1\times ...\times X^n=[0,1]^{[n]}=[0,1]^n.
\end{eqnarray*} 
The points $x=(x_1,...,x_n)\in X$ are the {\it mixed 
strategy combinations}.
The map $V^i:X\to\R^n$ is the multilinear extension of $U^i$.
It is a polynomial in $\R[x_1,...,x_n]$ where each variable
$x_i$ has in each monomial degree 0 or 1.
Then $V=(V^1,...,V^n):X\to\R^n$. 

We also denote 
$X^{-i}:=X^1\times ...\times X^{i-1}\times X^{i+1}\times ...
\times X^m$ and its elements
$x_{-i}:=(x_1,...,x_{i-1},x_{i+1},...,x_m)\in X^{-i}$.
We follow the standard (slightly incorrect)
convention and identify $X^i\times X^{-i}$ with $X$
and $(x_i,x_{-i})$ with $x$. 

\medskip
(c) For $([n],X,V)$ as above, define a function
$\Delta^i:X^{-i}\to\R$ for $i\in[n]$ by 
\begin{eqnarray}\label{2.1}
\Delta^i(x_{-i}):=V^i(1,x_{-i})-V^i(0,x_{-i})
\quad\textup{for }x_{-i}\in X^{-i}.
\end{eqnarray}
Then $x\in X$ is a Nash equilibrium if and only if 
for any $i\in[n]$
\begin{eqnarray}\label{2.2}
\Delta^i(x_{-i})=0&&
\textup{if }x_i\in (0,1),\\
\Delta^i(x_{-i})\geq 0&&
\textup{if }x_i=1,\label{2.3}\\
\textup{and }\Delta^i(x_{-i})\leq 0 &&
\textup{if }x_i= 0. \label{2.4}
\end{eqnarray}
The set $\NN$ is the set of all Nash equilibria. 
A {\it TMNE} ({\it totally mixed Nash equilibrium})
is a Nash equilibrium $x$ in the interior of $X$.

\medskip
(d) For $([n],X,V)$ as above, a mixed strategy combination
$x\in X$ is an {\it equilibrium candidate} if it satisfies
\eqref{2.2} (but not necessarily \eqref{2.3} and \eqref{2.4}).
\end{definition}

\begin{remarks}\label{t2.2}
(i) In Definition \ref{t2.1} (c) $\Delta^i$ is a polynomial
in $\R[x_1,...,x_{i-1},x_{i+1},...,x_n]$ where each variable
$x_j$ (with $j\neq i$) has in each monomial degree 0 or 1.

(ii) Vice versa, for any tuple $(\Delta^1,...,\Delta^n)$
of such polynomials a two-action game exists such that 
$\Delta^1,...,\Delta^n$ satisfy \eqref{2.1}.

(iii) In fact, one can choose this game such that
$V^i(0,x_{-i})=0$ for any $x_{-i}\in X^{-i}=[0,1]^{[n]-\{i\}}$.
Then $V^i(x)=x_i\cdot \Delta^i(x_{-i})$. 

(iv) The conditions in \eqref{2.3} and \eqref{2.4} can be put
together as
\begin{eqnarray}\label{2.5} 
(-1)^{1+x_i}\Delta^i(x_{-i})\geq 0\quad\textup{for }
x_i\in\{0;1\}.
\end{eqnarray}

(v) In a two-action game with $n$ players, the set
$X=[0,1]^{[n]}$ is a hypercube of dimension $n$. 
For $l\in\{0,1,...,n\}$ it has 
$\begin{pmatrix}n\\l\end{pmatrix}\cdot 2^l$ 
faces of dimension $n-l$, and they are also hypercubes. 

(vi) In a two-action game with $n$ players, each equilibrium
candidate $x$ is in the interior of a unique face, and it is a TMNE
of the restriction of the game to this face. This fact, part (v) 
and the upper bound of \nocite{McKelvey and McLennan (1997)}
show that in a generic two-action game the number of equilibrium
candidates is bounded from above by
$V(n)=\sum_{l=0}^n\begin{pmatrix}n\\l\end{pmatrix}\cdot 2^l\cdot 
!(n-l)$. This implies the upper bound 
$\mu(n;2,...,2)\leq V(n)$ of \nocite{Vidunas (2017)} for the 
number of Nash equilibria in a generic two-action game.

(vii) Lemma \ref{t3.4} (b)(iii) will present special two-action 
games where the number of equilibrium candidates is equal to $V(n)$.
Also for generic two-action games close to these special games
the number of equilibrium candidates is equal to $V(n)$. 
Therefore the maximal number of equilibrium candidates within
all generic two-action games with $n$ players is equal to $V(n)$. 
\end{remarks}

Throughout the whole paper, permutations without fixed points
will be important. We recall the corresponding notions.

\begin{definition}\label{t2.3}
(a) A permutation $\pi\in S_n$ (for some $n\in\N$) is a
{\it derangement} if $\pi(i)\neq i$ for any $i\in[n]$.
The set of derangements in $S_n$ is called $\textup{Der}_n$. 
The {\it subfactorial} $!n$ is the number $|\textup{Der}_n|$
of derangements in $S_n$.

(b) For $\pi\in S_n$, denote by $\Fix(\pi):=\{i\in [n]\,|\,
\pi(i)=i\}\subset[n]$ the set of fixed points of $\pi$. 
Then $\pi$ restricted to $[n]-\Fix(\pi)$ is a derangement 
on the set $[n]-\Fix(\pi)$.
\end{definition}

The subfactorial $!n$ can be determined in different ways.
The following lemma collects three ways. It is completely
elementary and well known, see e.g. Proposition 5.4 in 
(\nocite{McKelvey and McLennan (1997)}) or Example 2.2.1
in (\nocite{Stanley (2012)})

\begin{lemma}\label{t2.4}
The subfactorials are determined by each of the recursions
\eqref{2.6} and \eqref{2.7} and by the closed formula
\eqref{2.8}. The table \eqref{2.9}
gives the first subfactorials and (for comparison) the
first factorials.
\begin{eqnarray}\label{2.6}
!0=1,\ !n&=&n\cdot !(n-1) + (-1)^n\quad\textup{for }n\geq 1,\\
!0=1,\ !1=0,\  
!n&=&(n-1)(!(n-1)+!(n-2))\quad\textup{for }n\geq 2,\hspace*{1cm}
\label{2.7}\\
!n&=&n!\cdot\sum_{j=0}^n\frac{(-1)^j}{j!}\Bigl(
=n!\cdot\sum_{j=2}^n\frac{(-1)^j}{j!}\Bigr).\label{2.8}
\end{eqnarray}
\begin{eqnarray}
\begin{array}{r|r|r|r|r|r|r|r}
n&0&1&2&3&4&5&6\\ \hline
!n&1&0&1&2&9&44&265\\ \hline
n!&1&1&2&6&24&120&720
\end{array}\label{2.9}
\end{eqnarray}
The closed formula \eqref{2.8} 
shows $\lim_{n\to\infty}\frac{!n}{n!}=e^{-1}$. 
\end{lemma}

\section{Product two-action games}\label{c3}
\setcounter{equation}{0}

The main point in this section is the study
of a special class of games, which we call 
{\it product two-action games}, where we have the maximal number 
of equilibrium candidates and where we can understand rather well
which of them are Nash equilibria.

\begin{definition}\label{t3.1}
(a) A two-action game $([n],X,V)$ is a 
{\it product two-action game},
if for each player $i\in[n]$ the difference function 
$\Delta^i(x_{-i})=V^i(1,x_{-i})-V^i(0,x_{-i})$ in \eqref{2.1}
has the shape
\begin{eqnarray}\label{3.1}  
\Delta^i(x_{-i})&=&\prod_{j\in[n]-\{i\}}(x_j-a^i_j)
\end{eqnarray}
with suitable $a^i_j\in (0,1)$ where $a^{i_1}_j\neq a^{i_2}_j$ 
for $j\in[n]$ and $i_1,i_2\in[n]-\{j\}$ with $i_1\neq i_2$ is
demanded. Such games exist by the Remarks \ref{t2.2} (ii) and (iii).
By Remark \ref{t2.2} (iii) one can even choose the utility function
of player $i$ as $V^i(x)=x_i\Delta^i(x_{-i})$. 

\medskip
(b) Let $([n],X,V)$ be a product two-action game.
For $j\in[n]$, the numbers $a^1_j,a^2_j,...,a^n_j$ can be ordered
by size, starting with the smallest number and ending with the
biggest number. Let $p_j(i)\in[n-1]$
be the position of $a^i_j$ in this ordering.
Then $p_j$ is a bijection from the set $[n]-\{j\}$ to the
set $[n-1]:=\{1,...,n-1\}$. 
Its definition can be rewritten as follows,  
\begin{eqnarray}\label{3.2}
0<a^{(p_j)^{-1}(1)}_j
<a^{(p_j)^{-1}(2)}_{j}
<...<a^{(p_j)^{-1}(n-1)}_{j}<1
\end{eqnarray}
$p_j$ is called {\it $j$-th associated order map} (of the 
product two-action game). 
\end{definition}

\begin{examples}\label{t3.2}
(i) In the case $n=4$, the four difference functions are
\begin{eqnarray*}
\Delta^1(x_{-1})=\Delta^1(x_2,x_3,x_4)&=& 
\hspace*{1.6cm}(x_2-a^1_2)(x_3-a^1_3)(x_4-a^1_4),\\
\Delta^2(x_{-2})=\Delta^2(x_1,x_3,x_4)&=& 
(x_1-a^2_1)\hspace*{1.6cm}(x_3-a^2_3)(x_4-a^2_4),\\
\Delta^3(x_{-3})=\Delta^3(x_1,x_2,x_4)&=& 
(x_1-a^3_1)(x_2-a^3_2)\hspace*{1.6cm}(x_4-a^3_4),\\
\Delta^4(x_{-4})=\Delta^4(x_1,x_2,x_3)&=& 
(x_1-a^4_1)(x_2-a^4_2)(x_3-a^4_3)\hspace*{1.6cm}.
\end{eqnarray*}

(ii) The product two-action games with fixed $n$ and with a 
fixed tuple $P:=(p_1,p_2,...,p_n)$ of associated order maps
have all the same properties. We could restrict to one
tuple $(a^i_j)_{i\in[n],j\in[n]-\{i\}}$ of values.
For example the tuple of values
$$a^i_j:=\frac{p_j(i)}{n}$$
works. In the case $n=4$ this gives 
\begin{eqnarray*}
\Delta^1(x_{-1})&=& \hspace*{2.1cm}
(x_2-\frac{p_2(1)}{4})(x_3-\frac{p_3(1)}{4})(x_4-\frac{p_4(1)}{4}),\\
\Delta^2(x_{-2})&=& (x_1-\frac{p_1(2)}{4})
\hspace*{2.1cm}(x_3-\frac{p_3(2)}{4})(x_4-\frac{p_4(2)}{4}),\\
\Delta^3(x_{-3})&=& (x_1-\frac{p_1(3)}{4})
(x_2-\frac{p_2(3)}{4})\hspace*{2.1cm}(x_4-\frac{p_4(3)}{4}),\\
\Delta^4(x_{-4})&=& (x_1-\frac{p_1(4)}{4})
(x_2-\frac{p_2(4)}{4})(x_3-\frac{p_3(4)}{4})\hspace*{2.1cm}.
\end{eqnarray*}
\end{examples}

\begin{remarks}\label{t3.3}
\nocite{McKelvey and McLennan (1997)} 
considered in chapter 4 
special finite games for arbitrary $n$ and $S$ which have
the maximal number of TMNE. Our product two-action games
are those games in (\nocite{McKelvey and McLennan (1997)}, ch.~4)
which are also two-action. 
For them the maximal number of TMNE is $!n$.
We will recover this special case of the result in 
(\nocite{McKelvey and McLennan (1997)}, ch.~4)
in Corollary \ref{t3.6} (e).
But the main point of Lemma \ref{t3.4}, Lemma \ref{t3.5} and
Corollary \ref{t3.6} is a simultaneous
control of {\it all} equilibrium candidates and of the question 
which of them are Nash equilibria.
\end{remarks}

First we determine the equilibrium candidates of a product
two-action game. Their number turns out to be the maximal number
$V(n)$ of equilibrium candidates for generic two-action games,
which was determined in Remark \ref{t2.2} (vi).

\begin{lemma}\label{t3.4}
Let $([n],X,V)$ be a product two-action game
with numbers $a^i_j\in (0,1)$ as in Definition \ref{t3.1}
and tuple $P=(p_1,...,p_n)$ of associated order maps.

(a) A mixed strategy combination $x\in X=[0,1]^{[n]}$ is an
equilibrium candidate if and only if there exists a 
permutation $\pi\in S_n$ such that
\begin{eqnarray}
x_j=a^{\pi(j)}_j&&\textup{for }j\in [n]-\Fix(\pi),\label{3.3}\\
x_j\in\{0;1\}&&\textup{for }j\in \Fix(\pi).\label{3.4}
\end{eqnarray}
This $x$ is determined by the permutation $\pi$ and by the set 
$L_0(x):=\{j\in[n]\,|\, x_j=0\}\subset\Fix(\pi)$.

(b) Let $EC(\pi)$ be the set of equilibrium candidates
with permutation $\pi\in S_n$ (i.e. \eqref{3.3} and \eqref{3.4}
hold for $x\in EC(\pi)$ with $\pi$).

(i) Then $|EC(\pi)|=2^{|\Fix(\pi)|}$. 

(ii) The set of all equilibrium candidates is 
$\bigcup_{\pi\in S_n}EC(\pi)$.

(iii) The number of all equilibrium candidates is $V(n)$ in \eqref{1.2}.
\end{lemma}

{\bf Proof:} (a) 
$\Rightarrow$: Let $x\in X=[0,1]^{[n]}$ be an equilibrium 
candidate with $L(x):=\{i\in[n]\,|\, x_i\in\{0;1\}\}$. 
For each $i\in [n]-L(x)$ the equality \eqref{2.2} and 
the definition of the product two-action game show
$$0=\Delta^i(x_{-i})=\prod_{j\in [n]-\{i\}}(x_j-a^i_j),$$
so (at least) one of the factors in this product must be 0.
Therefore for each $i\in[n]-L(x)$ a value 
$\www{\pi}(i)\in [n]-\{i\}$ with 
$x_{\www{\pi}(i)}=a^i_{\www{\pi}(i)}$ exists. 

As all $a^i_j\in (0,1)$, $x_{\www{\pi}(i)}\in(0,1)$ and 
$\www{\pi}(i)\notin L(x)$.
As for fixed $j$ the $a^k_j$ for $k\in[n]-\{j\}$ are pairwise
different, $\www{\pi}(i_1)\neq\www{\pi}(i_2)$ for $i_1\neq i_2$.
Therefore $\www{\pi}$ is a derangement of the set $[n]-L(x)$.

It extends to a permutation $\www{\pi}\in S_n$ with
$\Fix(\www{\pi})=L(x)$ (by setting $\www{\pi}(i)=:i$ for
$i\in L(x)$). Let $\pi:=\www{\pi}^{-1}\in S_n$ be the inverse
permutation. Then also $\Fix(\pi)=L(x)$, and then 
\eqref{3.3} and \eqref{3.4} hold.

$\Leftarrow$: This is obvious.

Of course, $x$ is determined by $\pi$ and $L_0(x)$. 

(b) The parts (i) and (ii) are clear. 

(iii) For each $l\in\{0,1,...,n\}$ the number of subsets
$L\subset[n]$ with $|L|=l$ is $\begin{pmatrix}n\\l\end{pmatrix}$.
For each such subset $L$, the number of permutations 
$\pi\in S_n$ with $\Fix(\pi)=L$ is $!(n-l)$. Therefore
$$|\bigcup_{\pi\in S_n}EC(\pi)|=\sum_{l=0}^n
\begin{pmatrix}n\\l\end{pmatrix}\cdot 2^l\cdot !(n-l)=V(n).
\hspace*{1cm}\Box$$

\bigskip
Next we introduce a {\it sign tuple} for each equilibrium
candidate. Its shape will tell when the equilibrium candidate
is a Nash equilibrium. 
In the following, the sign function $\sign:\R^*\to\{\pm 1\}$
maps positive numbers to $+1$ and negative numbers to $-1$.

\begin{lemma}\label{t3.5}
Let $([n],X,V)$, $a^i_j$ and $P$ be as in Lemma \ref{t3.4}.

(a) Let $x\in X$ be an equilibrium candidate with permutation
$\pi\in S_n$. For $i\in\Fix(\pi)$ the number
$(-1)^{1+x_i}\Delta^i(x_{-i})$ in \eqref{2.5} is not zero.
Its sign $$s_i(x):=\sign((-1)^{1+x_i}\Delta^i(x_{-i}))$$
is called $i$-th sign of $x$.
The sign tuple of $x$ is 
$$s(x):=(s_i(x))_{i\in\Fix(\pi)}\in \{\pm 1\}^{\Fix(\pi)}.$$

(b) An equilibrium candidate $x$ is a Nash equilibrium
if and only if its sign tuple consists only of $+1$'s.

(c) Let $x\in X$ be an equilibrium candidate with
permutation $\pi$ with $\Fix(\pi)\neq\emptyset$,
and let $i\in\Fix(\pi)$. Then
\begin{eqnarray}
s_i(x)&=&(-1)^{|L_0(x)|}\cdot\prod_{j\in[n]-\Fix(\pi)}
\sign(a^{\pi(j)}_j-a^i_j)\nonumber \\
&=& (-1)^{|L_0(x)|}\cdot\prod_{j\in [n]-\Fix(\pi)}
\sign(p_j(\pi(j))-p_j(i)).\label{3.5}
\end{eqnarray}
\end{lemma}

{\bf Proof:} (a) The equilibrium candidate $x$ satisfies
\eqref{3.3} and \eqref{3.4}.  
Write $L_1(x):=\{i\in[n]\,|\, x_i=1\}=\Fix(\pi)-L_0(x)$. 
For $i\in Fix(\pi)$
\begin{eqnarray*}
\Delta^i(x_{-i})
&=&\prod_{j\in [n]-\{i\}}(x_j-a^i_j)\\
&=&\prod_{j\in L_0(x)-\{i\}}(-a^i_j)
\prod_{j\in L_1(x)-\{i\}}(1-a^i_j)
\prod_{j\in [n]-\Fix(\pi)}(a^{\pi(j)}_j-a^i_j)\\
&\neq&0.
\end{eqnarray*}

(b) This is true because \eqref{2.5}, namely 
$(-1)^{1+x_i}\Delta^i(x_{-i})\geq 0$ for $x_i\in\{0;1\}$, 
is equivalent to \eqref{2.3} and \eqref{2.4} by Remark 
\ref{t2.2} (iv). 

(c) The calculation in the proof of part (a) shows
for $i\in\Fix(\pi)$ 
\begin{eqnarray*}
s_i(x)=(-1)^{1+x_i}(-1)^{|L_0(x)-\{i\}|}
\prod_{j\in[n]-\Fix(\pi)}\sign(a^{\pi(j)}_j-a^i_j).
\end{eqnarray*}
The first two factors combine to $(-1)^{|L_0(x)|}$ in any case 
(whether or not $i\in L_0(x)$). The product of the other factors 
is because of the definition of the 
order map $p_j$ for the numbers $a^k_j$ with 
fixed $j$ and $k\in[n]-\{j\}$ equal to the product 
$\prod_{j\in[n]-\Fix(\pi)} \sign(p_j(\pi(j))-p_j(i))$. \hfill$\Box$

\bigskip
The following corollary puts together results on the 
Nash equilibria in $EC(\pi)$ for each $\pi\in S_n$.

\begin{corollary}\label{t3.6}
Let $([n],X,V)$, $a^i_j$ and $P$ be as in Lemma \ref{t3.4}.

(a) Consider $\pi\in S_n$. The element $x^{(\pi)}\in X$ with
$x^{(\pi)}_i=1$ for $i\in\Fix(\pi)$ and 
$x^{(\pi)}_i=a^{\pi(i)}_i$ for $i\in[n]-\Fix(\pi)$ is in
$EC(\pi)$. 

(b) Consider $\pi\in \Der_n$. The set $EC(\pi)$ contains only
$x^{(\pi)}$, and this is a Nash equilibrium.

(c) Consider $\pi\in S_n-\Der_n$. If $s(x^{(\pi)})$ is 
$(+1,...,+1)$ or $(-1,...,-1)$ then half of the elements of
$EC(\pi)$ are Nash equilibria. Else $EC(\pi)$ does not contain
Nash equilibria.

(d) (i) Consider $\pi\in S_n$ with $|\Fix(\pi)|=1$. 
Then $EC(\pi)$ has two elements, and one of them is a Nash
equilibrium.

(ii) There is no permutation $\pi\in S_n$ with $|\Fix(\pi)|=n-1$. 

(iii) Consider $\pi\in S_n$ with $|\Fix(\pi)|=n$. Then 
$\pi=\id$. The set $EC(\id)$ is the set of the $2^n$ vertices of
$X$, and half of its elements are Nash equilibria.

(e) The second line of the following table gives the number of 
Nash equilibria in the union of all open faces in $X$ of dimension 
$k\in\{n,n-1,1,0\}$.
It also gives an upper bound for this number for $k$ with
$2\leq k\leq n-2$.
\begin{eqnarray}\label{3.6}
\begin{array}{r|c|c|c|c|c}
k & n & n-1 & k\textup{ with }2\leq k\leq n-2 & 1 & 0 \\ \hline
\sharp\textup{Nash eq.} & !n & !(n-1)\cdot n & 
\leq \begin{pmatrix}n\\n-k\end{pmatrix}\cdot 2^{n-k-1}\cdot !k & 
0 & 2^{n-1}
\end{array}
\end{eqnarray}
\end{corollary}

{\bf Proof:}
(a) This follows from Lemma \ref{t3.4} (a).

(b) Of course $EC(\pi)$ contains only the element $x^{(\pi)}$. 
Its sign tuple is empty, so the condition that all its entries are 
$+1$ is satisfied emptyly. In other words, the conditions 
\eqref{2.3} and \eqref{2.4} are empty here, and \eqref{2.2} 
is satisfied.

(c) Let $x\in EC(\pi)$. Then for 
$j\in[n]-\Fix(\pi)$ $x_j=x_j^{(\pi)}=a^{\pi(j)}_j$. So $x$ and
$x^{(\pi)}$ can differ only in the indices $j\in \Fix(\pi)$. 
Formula \eqref{3.5} shows
\begin{eqnarray*}
s(x)=(-1)^{|L_0(x)|}\cdot s(x^{(\pi)}),
\end{eqnarray*}
so on $EC(\pi)$ only two opposite sign tuples arise,
and each arises for half of the elements of $EC(\pi)$. 
If one of them is $(+1,...,+1)$ then half of the elements
of $EC(\pi)$ are Nash equilibria. If none of them is 
$(+1,...,+1)$, then $EC(\pi)$ does not contain Nash equilibria.

(d) (i) If $s(x^{(\pi)})=+1$ then $x^{(\pi)}$ is a Nash equilibrium.
If $s(x^{(\pi)})=-1$ then the other element of $EC(\pi)$ is a 
Nash equilibrium.

(ii) This is obvious.

(iii) Obviously $\pi=\id$. Here $s(x^{(\pi)})=(+1,...,+1)$ by
\eqref{3.5}. One applies part (c).

(e) For a fixed $k\in\{0,1,...,n\}$ an equilibrium candidate on
an open face of dimension $k$ is in the set $EC(\pi)$ for some
$\pi\in S_n$ with $|\Fix(\pi)|=n-k$. Recall
$|EC(\pi)|=2^{n-k}$. Recall also that the number of permutations
$\sigma\in S_n$ with $|\Fix(\sigma)|=n-k$ is 
$\begin{pmatrix}n\\n-k\end{pmatrix}\cdot !k$. 
One applies part (b)
for $k=n$, part (d) for $k\in\{n-1,1,0\}$, and part (c) for $k$
with $2\leq k\leq n-2$. \hfill$\Box$

\bigskip

Finally, we consider a product two-action game of $n=4$ players
with a naively chosen tuple of associated orders.
We will see that not for all sets $EC(\pi)$ with $\pi\in S_4$
and $|\Fix(\pi)|=2$ half of their elements are Nash equilibria.

\begin{example}\label{t3.7} Consider a product two-action game
with $n=4$. Table \eqref{3.6} specializes to Table \eqref{3.7},
\begin{eqnarray}\label{3.7}
\begin{array}{r|c|c|c|c|c}
k & 4 & 3 & 2 & 1 & 0 \\ \hline
\textup{number of Nash equilibria} & 9 & 8 & \leq 12 & 0 & 8
\end{array}
\end{eqnarray}
For the faces of dimension 2, the number of Nash equilibria
on these faces depends on the values $a^i_j$. 
We choose $a^i_j=\frac{p_j(i)}{4}$, 
and we choose the order maps $p_1,p_2,p_3$ and $p_4$ as simple as
possible, namely as follows.
\begin{eqnarray*}
\begin{array}{c|c|c|c}
p_1 & p_2 & p_3 & p_4 \\
2\ 3\ 4 & 1\ 3\ 4 & 1\ 2\ 4 & 1\ 2\ 3 \\
\downarrow\ \downarrow\ \downarrow &
\downarrow\ \downarrow\ \downarrow &
\downarrow\ \downarrow\ \downarrow &
\downarrow\ \downarrow\ \downarrow \\
1\ 2\ 3 & 1\ 2\ 3 & 1\ 2\ 3 & 1\ 2\ 3
\end{array}
\end{eqnarray*}
There are six subsets $L\subset[4]$ with $|L|=2$, and for
each such subset $L$ there is one permutation $\pi\in S_4$
with $\Fix(\pi)=L$. It is the transposition of the two 
elements in $[4]-\Fix(\pi)$. 
In each of these six cases we consider
$x^{(\pi)}\in EC(\pi)$ with $x_j^{(\pi)}=1$ for $j\in L$ and 
$x^{(\pi)}_j=a^{\pi(j)}_j=\frac{p_j(\pi(j))}{4}$ for $j\in[4]-L$.
We also calculate its sign tuple 
$s(x^{(\pi)})=(s_i(x^{(\pi)})_{i\in\Fix(\pi)}$. The entries are
\begin{eqnarray*}
s_i(x^{(\pi)})&\stackrel{\textup{\eqref{3.5}}}{=}
&\prod_{j\in [4]-L}\sign(p_j(\pi(j))-p_j(i))\\
&\stackrel{\textup{here}}{=}&\prod_{j\in [4]-L}\sign(\pi(j)-i)
=\prod_{j\in [4]-L}\sign(j-i).
\end{eqnarray*}
\begin{eqnarray*}
\begin{array}{cccl}
L & \pi & x^{(\pi)} & s(x^{(\pi)}) \\ \hline 
\{1,2\} & (3\ 4) & (1,1,\frac{3}{4},\frac{3}{4}) & 
(s_1,s_2)=(1,1) \\
\{1,3\} & (2\ 4) & (1,\frac{3}{4},1,\frac{2}{4}) & 
(s_1,s_3)=(1,-1) \\
\{1,4\} & (2\ 3) & (1,\frac{2}{4},\frac{2}{4},1) & 
(s_1,s_4)=(1,1) \\
\{2,3\} & (1\ 4) & (\frac{3}{4},1,1,\frac{1}{4}) & 
(s_2,s_3)=(-1,-1) \\
\{2,4\} & (1\ 3) & (\frac{2}{4},1,\frac{1}{4},1) & 
(s_2,s_4)=(-1,1) \\
\{3,4\} & (1\ 2) & (\frac{1}{4},\frac{1}{4},1,1) & 
(s_3,s_4)=(1,1)
\end{array}
\end{eqnarray*}
Therefore $EC(\pi)$ for $\pi\in\{(3\ 4),(2\ 3),(1\ 4),
(1\ 2)\}$ contains two Nash equilibria, but $EC(\pi)$
for $\pi\in\{(2\ 4),(1\ 3)\}$ contains no Nash equilibria.
Section \ref{c4} will provide a different
tuple of associated order maps such that then $EC(\pi)$
for each transposition $\pi$ contains two Nash equilibria.
See Example \ref{t4.3}.
\end{example}

\section{Maximal product two-action games}\label{c4}
\setcounter{equation}{0}

This section proves the main result of the paper,
the existence of generic two-action games
with $\frac{1}{2}(V(n)+!n)$ Nash equilibria. 
The main point is to find product two-action games with
so many Nash equilibria.

Recall that Corollary \ref{t3.6} (c) implies for 
a product two-action game $([n],X,V)$ and a permutation
$\pi\in S_n-\Der_n$ the following:
Either half of the equilibrium candidates in $EC(\pi)$ are
equilibria or none of them are equilibria
(none is impossible in the cases $\pi=\id$ and $|\Fix(\pi)|=1$
by Corollary \ref{t3.6} (d)).

\begin{definition}\label{t4.1}
A product two-action game $([n],G,V)$ is {\sf maximal}
if for each permutation $\pi\in S_n-\textup{Der}_n$ 
half of the equilibrium candidates in $EC(\pi)$ are equilibria.
\end{definition}

The main result of this paper is that for any number 
$n\in\N$ maximal product two-action games exists.
The point is to find a suitable tuple $(p_1,...,p_n)$
of associated order maps. Theorem \ref{t4.2} says that the following
tuple $(p_1^*,...,p_n^*)$ of associated order maps does it.
Working out Lemma \ref{t3.4}, Lemma \ref{t3.5} and
Corollary \ref{t3.6} was straightforward. In contrast, 
finding the tuple $(p_1^*,...,p_n^*)$ in Theorem \ref{t4.2}
was not straightforward, but required a lot of effort.

\begin{theorem}\label{t4.2}
Let $n\in\Z_{\geq 2}$. A product two-action game with $n$ players
and with the following tuple $(p_1^*,...,p_n^*)$ of associated 
order maps is maximal. For $j\in[n]$ the bijection 
$p_j^*:[n]-\{j\}\to [n-1]$ is defined by 
\begin{eqnarray*}
p_j^*:\left(\begin{array}{ll}
i\mapsto i+n-j& \textup{ if }1\leq i\leq j-1\\
i\mapsto i-j  & \textup{ if }j+1\leq i\leq n
\end{array}\right).
\end{eqnarray*}
\end{theorem}

The following picture illustrates $p_j^*$:

\centerline{
\begin{tikzpicture}[scale=0.8]
\draw [thick] (0,3)--(4,3);
\draw [thick] (0,3)--(0,3.3);
\draw [thick] (4,3)--(4,3.3);
\node at (0.2,3.3) {1};
\node at (3.7,3.3) {j-1};
\node at (2,3.3) {$\cdots$};
\draw [-stealth] (2,2.9)--(8.2,2.1);
\draw [thick] (5,3)--(11,3);
\draw [thick] (5,3)--(5,3.3);
\draw [thick] (11,3)--(11,3.3);
\node at (5.5,3.3) {j+1};
\node at (10.7,3.3) {n};
\node at (8,3.3) {$\cdots$};
\draw [-stealth] (8,2.9)--(3,2.1);
\draw [thick] (0,1.5)--(6,1.5);
\draw [thick] (0,1.5)--(0,1.8);
\draw [thick] (6,1.5)--(6,1.8);
\node at (0.2,1.8) {1};
\node at (5.6,1.8) {n-j};
\node at (3,1.8) {$\cdots$};
\draw [thick] (6.2,1.5)--(10.2,1.5);
\draw [thick] (6.2,1.5)--(6.2,1.8);
\draw [thick] (10.2,1.5)--(10.2,1.8);
\node at (6.9,1.8) {n-j+1};
\node at (9.8,1.8) {n-1};
\node at (8.2,1.8) {$\cdots$};
\end{tikzpicture}
}

{\bf Proof:}
Consider a product two-action game with the tuple
$(p_1^*,...,p_n^*)$ of associated order maps.
Fix a permutation $\pi\in S_n$ with $2\leq|\Fix(\pi)|\leq n-2$.
Table \eqref{t3.6} in Corollary \ref{t3.6} covers the other cases.

Consider the equilibrium candidate $x^{(\pi)}\in EC(\pi)$ 
with $L_0(x^{(\pi)})=\emptyset$, so with $x_i=1$ for each
$i\in \Fix(\pi)$. By Corollary \ref{t3.6} (c),
it is sufficient to show that the sign tuple
$s(x^{(\pi)})$ is either $(+1,...,+1)$ or $(-1,...,-1)$.
For $i\in \Fix(\pi)$ the sign $s_i(x^{(\pi)})$ is 
\begin{eqnarray*}
s_i(x^{(\pi)})\stackrel{\textup{\eqref{3.5}}}{=}
\prod_{j\in[n]-\Fix(\pi)}\sign(p_j^*(\pi(j))-p_j^*(i)).
\end{eqnarray*}

The following table gives for $i\in\Fix(\pi)$ and 
$j\in[n]-\Fix(\pi)$ for all possible orders of the values
$i,j$ and $\pi(j)$ the sign $\sign(p_j^*(\pi(j))-p_j^*(i))$.
\begin{eqnarray}\label{4.1}
\begin{array}{ccccc|l}
&&&& & \sign(p_j^*(\pi(j))-p_j^*(i))\\ \hline
\pi(j)<&j&<&i& & +1 \\
&j&<\pi(j)<&i& & -1 \\
&j&<&i&<\pi(j) & +1 \\
\pi(j)<&i&<&j& & -1 \\
&i&<\pi(j)<&j& & +1 \\
&i&<&j&<\pi(j) & -1 
\end{array}
\end{eqnarray}
The table says
\begin{eqnarray*}
s_i(x^{(\pi)})=(-1)^{|\{j\,|\, j<\pi(j)<i\}|}
(-1)^{|\{j\,|\, \pi(j)<i<j\}|}
(-1)^{|\{j\,|\, i<j<\pi(j)\}|}.
\end{eqnarray*}
Recall $i\in\Fix(\pi)$. Of course
\begin{eqnarray*}
|\{j\in[n]\,|\, j<i<\pi(j)\}| 
=|\{j\in[n]\,|\, \pi(j)<i<j\}|.
\end{eqnarray*}
Therefore
\begin{eqnarray*}
s_i(x^{(\pi)})&=&(-1)^{|\{j\,|\, j<\pi(j)<i\}|}
(-1)^{|\{j\,|\, j<i<\pi(j)\}|}
(-1)^{|\{j\,|\, i<j<\pi(j)\}|}\\
&=& (-1)^{|\{j\,|\, j<\pi(j)\}|}
\end{eqnarray*}
This is independent of $i$. Therefore the sign tuple of $x^{(\pi)}$
is $s(x^{(\pi)})=(-1)^{|\{j\,|\, j<\pi(j)\}|}(+1,...,+1)$, 
so it is either
$(+1,...,+1)$ or $(-1,...,-1)$. This finishes the proof of
Theorem \ref{t4.2}. \hfill$\Box$

\begin{example}\label{t4.3} 
Consider the product two-action game with $n=4$, with tuple
$(p_1^*,p_2^*,p_3^*,p_4^*)$ of associated order maps
and with $a^i_j=\frac{p_j^*(i)}{4}$. 
These order maps are as follows.
\begin{eqnarray*}
\begin{array}{c|c|c|c}
p_1^* & p_2^* & p_3^* & p_4^* \\
2\ 3\ 4 & 1\ 3\ 4 & 1\ 2\ 4 & 1\ 2\ 3 \\
\downarrow\ \downarrow\ \downarrow &
\downarrow\ \downarrow\ \downarrow &
\downarrow\ \downarrow\ \downarrow &
\downarrow\ \downarrow\ \downarrow \\
1\ 2\ 3 & 3\ 1\ 2 & 2\ 3\ 1 & 1\ 2\ 3
\end{array}
\end{eqnarray*}
The following table gives for each subset $L\subset[4]$ with 
$|L|=2$ and for the transposition $\pi\in S_4$ with $\Fix(\pi)=L$
the equilibrium candidate $x^{(\pi)}\in EC(\pi)$ and its sign tuple. 
Then $x^{(\pi)}_j=1$ for $j\in L$ and 
$x^{(\pi)}_j=a^{\pi(j)}_j=\frac{p_j^*(\pi(j))}{4}$ for 
$j\in [4]-L$. By the proof of Theorem \ref{t4.2}
$s_i(x^{(\pi)})=(-1)^{|\{j\,|\, j<\pi(j)\}|}=-1$ for $i\in \Fix(\pi)$. 
One can check this directly by applying formula \eqref{3.5},  
\begin{eqnarray*}
s_i(x^{(\pi)})&\stackrel{\textup{\eqref{3.5}}}{=}
&\prod_{j\in [4]-L}\sign(p_j^*(\pi(j))-p_j^*(i)).
\end{eqnarray*}

\begin{eqnarray*}
\begin{array}{cccl}
L & \pi & x^{(\pi)} & s(x^{(\pi)}) \\ \hline 
\{1,2\} & (3\ 4) & (1,1,\frac{1}{4},\frac{3}{4}) & 
(s_1,s_2)=(-1,-1) \\
\{1,3\} & (2\ 4) & (1,\frac{2}{4},1,\frac{2}{4}) & 
(s_1,s_3)=(-1,-1) \\
\{1,4\} & (2\ 3) & (1,\frac{1}{4},\frac{3}{4},1) & 
(s_1,s_4)=(-1,-1) \\
\{2,3\} & (1\ 4) & (\frac{3}{4},1,1,\frac{1}{4}) & 
(s_2,s_3)=(-1,-1) \\
\{2,4\} & (1\ 3) & (\frac{2}{4},1,\frac{2}{4},1) & 
(s_2,s_4)=(-1,-1) \\
\{3,4\} & (1\ 2) & (\frac{1}{4},\frac{3}{4},1,1) & 
(s_3,s_4)=(-1,-1)
\end{array}
\end{eqnarray*}
Only the sign tuple $(-1,-1)$ arises.
Therefore for each transposition $\pi$, the set $EC(\pi)$ 
contains two Nash equilibria. The product two-action game
is maximal. 
\end{example}

\begin{corollary}\label{t4.4}
(a) A maximal product two-action game with $n$ players has
$\frac{1}{2}(V(n)+!n)$ Nash equilibria.

(b) A small generic deformation of a product two-action
game has the same number of Nash equilibria as the
product two-action game.

(c) The inequality \eqref{1.3}
$\frac{1}{2}(V(n)+!n)\leq \mu(n;2,...,2)$ holds.
\end{corollary}

{\bf Proof:}
(a) In a maximal product two-action game, the number of Nash
equilibria is
\begin{eqnarray*}
\sum_{\pi\in \Der_n}1
+\sum_{\pi\in S_n-\Der_n}\frac{1}{2}|EC(\pi)|
=\frac{1}{2}!n+\frac{1}{2}\sum_{\pi\in S_n}|EC(\pi)|
=\frac{1}{2}(!n+V(n)).
\end{eqnarray*}
(b) An equilibrium candidate $x\in EC(\pi)$ 
of a product two-action game satisfies 
$L_0(x)\dot\cup L_1(x)=\Fix(\pi)$, and it is an intersection 
point of the zero hypersurfaces of the functions in the following
tuple of functions,
\begin{eqnarray}\label{4.2}
\Bigl(x_i\,|\, i\in L_0(x);\quad
x_i-1\,|\, i\in L_1(x); \quad
\Delta^i(x_{-i})\,|\,
i\in [n]-\Fix(\pi)\Bigr).
\end{eqnarray}
Locally near $x$, these zero hypersurfaces coincide with
the zero hypersurfaces of the functions in the following tuple,
\begin{eqnarray}\label{4.3}
\Bigl(x_i\,|\, i\in L_0(x);\quad
x_i-1\,|\, i\in L_1(x); \quad
x_i-a_i^{\pi(i)}\,|\,
i\in [n]-\Fix(\pi)\Bigr).
\end{eqnarray}
Obviously locally near $x$, these zero hypersurfaces are smooth 
and transversal (therefore $x$ may be called a {\it regular} 
equilibrium candidate, see e.g. (\nocite{Harsanyi (1973)}). 
This property is preserved for the zero hypersurfaces
in the tuple in \eqref{4.2} by a small deformation
of the product two-action game to a generic two-action game.
Therefore $x$ deforms to an equilibrium candidate
of the deformed game. 

If $x$ is a Nash equilibrium, it satisfies additionally
the inequalities
\begin{eqnarray*}
\Delta^i(x_{-i})&<&0\quad\textup{for }i\in L_0(x),\\
\Delta^i(x_{-i})&>&0\quad\textup{for }i\in L_1(x).
\end{eqnarray*}
Also these inequalities are preserved by a small generic deformation.
Therefore a Nash equilibrium deforms to a Nash equilibrium.

(c) One considers a small generic deformation of a maximal
product two-action game and applies (a) and (b).
\hfill$\Box$

\bigskip
We conjecture (optimistically) that the maximal number of
Nash equilibria in product two-action games is also the
maximal number in all generic two-action games.

\begin{conjecture}\label{t4.5}
$\mu(n;2,...,2)=\frac{1}{2}(V(n)+!n)$ for $n\geq 2$. 
\end{conjecture}

This is well known to be true in the case $n=2$.
In the case $n=3$ it was recently proved
by Jahani and von Stengel (2024): $\mu(3;2,2,2)=9$.

\begin{remarks}\label{t4.6}
(i) The equality $\mu(3;2,2,2)=9$ is also 
claimed by \nocite{Vidunas (2017)} in Remark 5.6, 
but without proof. 

(ii) The case of a two-action game with $n=3$ is also
considered in the following three papers. 

\nocite{Chin, Parthasarathy and Raghavan (1973)} 
prove in Theorem 6 for each two-action game with $n=3$ and with
$\NN\subset\textup{interior}(X)$ that $|\NN|=1$.
For generic games this follows from the upper bound $2$
for the number of TMNE and from the oddness of $|\NN|$.
But they proved it for each two-action game with $n=3$.

\nocite{McKelvey and McLennan (1997)} in chapter 6 
reprove in an explicit way their
upper bound $!3=2$ for the number of TMNE in a generic
two-action game with $n=3$.

\nocite{Jahani and von Stengel (2022)} present an 
algorithmic approach to finding all Nash equilibria of any 
two-action game with $n=3$. 
\end{remarks}

\end{document}